# An exact relaxation of AC-OPF problem for battery-integrated power grids

H. Sekhavatmanesh, *Member, IEEE* and S. Mastellone, *Member, IEEE.*

*Abstract*—Renewable energy resources and power electronics-interfaced loads introduce fast dynamics in distribution networks. These dynamics cannot be regulated by slow conventional solutions and require fast controllable energy resources such as Battery Energy Storage Systems (BESSs). To compensate for the high costs associated to BESSs, their energy and power management should be optimized. In this paper, a convex iterative optimization approach is developed to find the optimal active and reactive power setpoints of BESSs in active distribution networks. The objective is to minimize the total cost of energy purchase from the grid. Round-trip and life-time characteristics of BESSs are modelled accurately and integrated into a relaxed and exact formulation of the AC power flow, resulting into a Modified Augmented Relaxed Optimal Power Flow (MAROP) problem. The feasibility and optimality of the solution under the grid security limits and technical constraints of BESSs is proven analytically. A 32-bus IEEE test benchmark is used to illustrate the performance of the developed approach in comparison to the alternative approaches existing in the literature.

## I. Introduction

The penetration of Renewable Energy Sources (RESs) in power grids is continuously increasing due to the growth of the electricity demand and sustainability requirements. Conventional operation mechanisms in distribution networks cannot compensate the RES frequent power intermittency, non-predictability and non-dispatchability. BESSs with fast-charging and -discharging capabilities can facilitate the integration of RESs into power grids. Besides capacity firming of RESs, many other applications are defined for BESSs in the power grids, including peak-power shaving, network congestion relief, grid upgrade deferral, emission reduction, voltage regulation, frequency services, and loss reduction [1]. Despite all the recent development in large-scale BESSs, they are still more expensive than conventional operation mechanisms. To compensate for the increased investment cost and realize all the BESS potential in the operation of power grids, it is essential to optimize the use of their power and energy capacities.

In this paper, we refer to the problem of finding the optimal active and reactive power setpoints of BESSs *as BESS operation planning problem*. For ensuring the solution feasibility, we require an accurate model of the BESS technical characteristics, including power density, energy density, lifetime, and round-trip efficiency. These characteristics together with those of the grid should be integrated into the optimization problem. In the following, we review the existing BESS and grid models.

The authors of [2] and [3] provide linear models of BESS, which they apply on the BESS operation planning problem and optimal home energy management, respectively. However, they assumed that BESS has perfect round-trip efficiencies and therefore do not model the non-negligible conversion losses of the BESS. The behavior of these losses is different in charging and discharging modes. To integrate them into the BESS model, different approaches are proposed in literature. In section V.B, we illustrate those models' shortcomings via simulation studies.

Another technical characteristic of the BESS is its lifetime, which limits the allowed number of discharging cycles. Inclusion of this constraint in the short-term BESS operation planning is rarely considered in the literature. As suggested by the authors of [4], in order to respect the maximum life-time cycle in long-term, a limit should be imposed on the number of discharging cycles in short-term operation of BESS. This limit should be a function of the desired end-of-life time of the BESS. This requires to determine the number of total discharging cycles during the planned time horizon as a function of the BESS power setpoints. This is formulated in [5], while disregarding round-trip efficiencies of the BESS, addressed then in [6], assuming that the energy level at the start and end of the day are equal to each other. This assumption is unrealistic for short-term planning of BESS. Moreover, all the BESS losses were not considered accurately.

Finally, the grid model should also be integrated into the BESS operation planning problem to ensure that the nodal voltage and line current magnitudes in the whole grid are bounded within the prescribed limits. To deal with the nonlinearity and nonconvexity inherent in the AC power flow equations, there are currently four approaches: I) Mixed-Integer Linear Optimization approaches, where the AC power flow formulations are either disregarded [7] or approximated in a linear way [8], II) rule-based approaches based on power flow simulations [9], III) non-linear optimization methodologies [10], IV) meta-heuristic approaches [11], and V) convexification methods [12]. The solutions provided by methods I, II, III, IV could be far from the global optimum. These methods are typically computationally heavy and might even fail to find an existing feasible solution. The convexification methods include relaxation and restriction methods that, respectively, provide lower and upper bounds of the objective function.

The relaxation methods, such as chordal relaxation [13], semidefinite programming (SDP) [14], [15], least squares estimation-based SDP [16], and *Second-Order Cone Programming* (SOCP) [17], [18] provide no guarantee for the solution feasibility. It is shown in [17], and [18] that the general SOCP relaxation method violates the upper voltage and current limits during high power production periods. In those papers, specific conditions are developed, under which the SOCP relaxation is exact and meets the security constraints. For example, the authors of [12] proved that if the objective function is strictly increasing with the line current magnitudes, then the solution of the SOCP problem is exact. This condition does not hold, e.g. in a case, where the objective is to track the reference of power injection from the upper grid [19].

The restriction methods provide solutions that are feasible but might be far from the global optimal one. According to [20]–[23], in these methods, the nonlinear power flow equations are linearized around a fixed operating point. In order to reduce the distance of these inner approximation envelopes from the exact equations, it is proposed in [21] to modify the fixed point through an iterative approach. However, no analytical proof is presented that guarantees the improvement of solution quality through iterations in case of each operation scenario.

The relaxation and restriction methods are combined in [24] to form Augmented Relaxed OPF (AROPF) formulation, ensuring the solution feasibility and optimality. In this approach, some conservative bounds on the voltage and line current magnitudes are added to the formulation of the SOCP relaxation method. An analytical proof is provided for feasibility and optimality of the solution under a priori verifiable conditions. However, the AROPF does not include the BESS model.

To formulate the BESS operation planning problem and overcome some of the modelling limitations, we propose a convex optimization problem. The objective function is to minimize the cost of power purchase from the grid while respecting all the BESS and grid security constraints. The paper contributions are as follows:

- Developing a convex model for the conversion losses and life-time constraints of the BESS that results in a non-conservative and feasible solution.
- For the integration of the proposed BESS model into the grid model, a convex iterative optimization algorithm is proposed.
- It is analytically proved that, under some conditions that can be checked a priori, the obtained solution of the optimization problem at each iteration satisfies all the grid and BESS constraints and the developed iterative approach converges to a non-conservative solution with respect to the modelling of the BESS conversion losses.

## II. Problem Statement

In this paper, we consider an Active Distribution Network, defined as an interconnected network of passive and active elements at the distribution level. Passive elements refer to noncontrollable elements, which include non-flexible loads and renewable energy sources such as Photovoltaic (PV) units. The active elements denote the units, whose active and reactive power setpoints can be regulated. The active elements that are considered in this paper include Dispatchable Distributed Generators (DDGs) (e.g. fuel cells and micro turbine units), and BESSs. All these elements are shown in a simple ADN in Fig. 1.

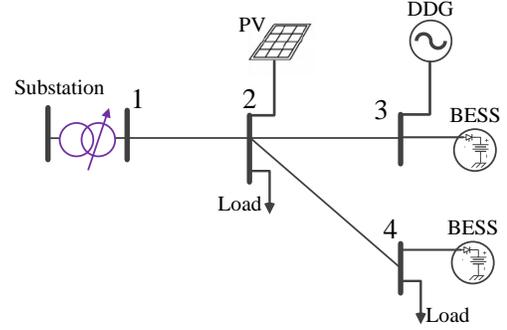

Fig. 1. A simple AND with active and passive elements considered in this paper.

As it can be seen, the ADN is connected to the transmission (or sub-transmission) grid through a substation bus. The voltage of this substation bus is fixed by the transmission system operator through a tap-changing transformer. Moreover, the network is assumed to be balanced, so a single-phase representation is used in Fig. 1 for the three-phase electrical network. In this paper, we define the optimal active and reactive power setpoints of active elements (i.e. DDGs and BESSs) in ADN so that the total power purchase from the grid in a short time horizon is minimized. To this end, a convex model for BESS operational planning problem is developed including all the constraints of the battery and grid.

## III. Modeling

### A. DGs' model

The RESs, and DDGs are regarded as negative loads and modelled as power-constant units. The power of loads and RESs are fixed and considered as parameters of the optimization problem, while the power of DDGs are considered as optimization variables. Together with the modeling of the BESS and grid, which are provided in the following of this section, the optimization problem can then be formulated.

### B. BESS' model

Fig. 2 shows the schematic of the BESS model considering its internal losses that include ohmic and conversion losses. We model maximum rate of power change, energy capacity, and power capacity, as BESS constraints. The non-linear original model of the BESS is formulated in (1).

$$E_{l,t} = E_{l,t=0} + \sum_{t^*=1}^{t} (p^r_{l,t^*} - r^b_l f^b_{l,t}) \Delta t \tag{1.a}$$

$$p^b_{l,t} = \begin{cases} \dfrac{p^r_{l,t}}{\eta^c_l}: & p^r_{l,t} \geq 0 \\ \eta^d_l p^r_{l,t}: & p^r_{l,t} \leq 0 \end{cases} = max\left\{\dfrac{p^r_{l,t}}{\eta^c}, \eta^d_l p^r_{l,t}\right\} \tag{1.b}$$

$$p^{b\,2}_{l,t} + q^{b\,2}_{l,t} = f^b_{l,t} v_{l,t} \tag{1.c}$$

$$E^{DSC}_l = \dfrac{\eta^c_l}{1 - \eta^c_l \eta^d_l} \sum_{t^*=1}^{T} \left(p^b_{l,t^*} - \dfrac{p^r_{l,t^*}}{\eta^c_l}\right) \tag{1.d}$$

where, $E_{l,t}$ is the State Of Charge (SOC) of the BESS at node $l$ at time $t$, $\Delta t$ is the optimization time step, $p^r_{l,t}$ is the active power absorbed to the reservoir of the BESS, regardless of the ohmic losses (see Fig. 2), $r^b_l$ is an equivalent resistance representing the internal power transmission losses of the BESS at node $l$, $f_{l,t}^b$ is the square of the current magnitude flowing into BESS at node $l$ and time $t$, $p_{l,t}^b$ and $q_{l,t}^b$ are active and reactive powers of the BESS at node $l$ that are exchanged with the grid at time $t$, $v_{l,t}$ is the square of the voltage magnitude at node $l$ and time $t$, $\eta_l^c$ and $\eta_l^d$ represent the charging and discharging power efficiencies of the BESS at node $l$, respectively, $E_l^{DSC}$ refers to the total discharged energy of the BESS at node $l$. $p_{l,t}^r$ and $p_{l,t}^b$ are positive if the BESS is in charging mode and negative if BESS is in discharging mode.

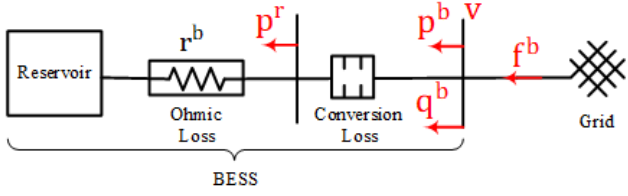

Fig. 2. The BESS' model.

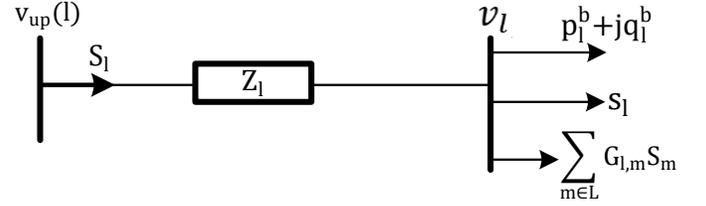

Fig. 3. The model of a distribution line with the notations used for the grid modelling

The SOC of a BESS at node $l$ and time $t$ is modelled in (1.a). The model of the conversion losses relates $p_{l,t}^b$ to $p_{l,t}^r$ as in (1.b). Equation (1.c) formulates the square of the current flow in the BESS. As proved in Appendix VII.D, (1.d) yields the total discharged energy of the BESS. This formulation is used later on to limit the number of the operation cycles for each BESS during the planning time horizon $T$ to a maximum limit, denoted as $N^{DSC}$. One charge cycle is completed when an amount of energy equal to 100% of the battery capacity is used, but not necessarily all from one charge. $N_c$ is chosen as a function of the desired end-of-life time of the BESS. For example, assume that a battery can perform 20000 cycles in its lifetime, and we aim to operate the BESS for 20 years. In such a case, $N^{DSC}$ for one day of operation is equal to 2.7.

The model of $p_{l,t}^b$ and $f_{l,t}^b$ in (1.b) and (1.c) are nonlinear and render the optimization problem non-convex. To address this issue, an upper and a lower bound are provided for $p_{l,t}^b$. The upper bound is denoted with $\bar{p}_{l,t}^b$ and formulated in (2.a).

$$\begin{cases} \bar{p}_{l,t}^b \geq \dfrac{p_{l,t}^r}{\eta_l^c} \\ \bar{p}_{l,t}^b \geq \eta_l^d p_{l,t}^r \end{cases} \quad (2.a)$$

$$\hat{p}_{l,t}^b = \begin{cases} \dfrac{p_{l,t}^r}{\eta_l^c}: & \tilde{p}_{l,t}^r \geq 0\,(charging\ state\ assumed) \\ \eta_l^d p_{l,t}^r: & \tilde{p}_{l,t}^r \leq 0\,(discharging\ state\ assumed) \end{cases} \quad (2.b)$$

$$\bar{p}_{l,t}^{b\,2} + q_{l,t}^{b\,2} \leq \bar{f}_{l,t}^b v_{l,t} \quad (2.c)$$

$$\underline{E}_{l,t} = E_{l,t=0} + \sum_{t^*=1}^{t}(p_{l,t^*}^r - r_l^b \bar{f}_{l,t}^b)\Delta t \quad (2.d)$$

$$\bar{E}_{l,t} = E_{l,t=0} + \sum_{t^*=1}^{t} p_{l,t^*}^r \Delta t \quad (2.e)$$

$$\bar{E}_l^{DSC} = \frac{\eta_l^c}{1 - \eta_l^c \eta_l^d} \sum_{t^*=1}^{T}\left(\bar{p}_{l,t^*}^b - \frac{p_{l,t^*}^r}{\eta_l^c}\right) \quad (2.f)$$

The formulation of $\hat{p}_{l,t}^b$ as the lower bound on $p_{l,t}^b$ is given in (2.b). Using $\tilde{p}_{l,t}^r$, we guess if at a given time $t$, the BESS is at charging or discharging state. To determine this guess, an iterative algorithm will be presented in section IV.D.

Lemma I: if $(p^r, p^b)$ satisfies (1.b) and $(p^r, \hat{p}^b)$ satisfies (2.b), we have $\hat{p}^b \leq p^b$ for any arbitrary choice of $\tilde{p}^r$.

The proof of Lemma I is provided in Appendix VII.A. The non-linear equation (1.c) is convexified introducing the auxiliary variable $\bar{f}_{l,t}^b$ as an upper bound on $f_{l,t}^b$, which is defined in (2.c) as a function of $\bar{p}_{l,t}^b$. The lower bound on $f_{l,t}^b$ is considered zero. Using the introduced upper and lower bounds on $f_{l,t}^b$, lower and upper-bounds on the energy level of the BESS expressed originally in (1.a) are developed in (2.d) and (2.e), respectively. Finally, (2.f) introduces $\bar{E}_l^{DSC}$ as an upper bound on the total discharged energy ($E_l^{DSC}$) expressed in (2.f), replacing $p_{l,t}^b$ with its upper bound $\bar{p}_{l,t}^b$.

### C. Grid's model

In this section, the electrical security constraints are formulated as functions of decision variables. In this regard, AC power flow equations are introduced in (3) according to the branch flow model.

The used notations are depicted in Fig. 3 for a typical distribution line. Index $l$ refers to a given line and also to its ending node. The starting node of line $l$ is indexed with $up(l)$. Index 0 corresponds to the slack bus. $s_{l,t}$ represents the net complex power consumption of the load, RES, and DDGs at node $l$. As shown in Fig. 3, the net power absorption at node $l$ is equal to $s_l + p_l^b + jq_l^b$. $S_{l,t}$ represents the apparent line power flow. The resistance, reactance, and impedance of line $l$ are expressed with $r_l$, $x_l$, and $z_l$, respectively. $\Re(.)$ and $\Im(.)$ denote the real and imaginary parts of a complex number, respectively.

$$S_{l,t} = s_{l,t} + p_{l,t}^b + jq_{l,t}^b + \sum_{m\in L} G_{l,m} S_{m,t} + Z_l f_{l,t} \quad (3.a)$$

$$v_{l,t} = v_{up(l),t} - 2\Re(z_l^* S_{l,t}) + |z_l|^2 f_{l,t} \quad (3.b)$$

$$\Re(S_{l,t})^2 + \Im(S_{l,t})^2 = f_{l,t} v_{up(l),t} \quad (3.c)$$

Equations (3.a) and (3.b) apply, respectively, current and voltage Kirchhoff's laws to a given line $l$ as shown in Fig. 3. Inequality (4) relaxes the non-linear equation given in (3.c). This relaxation might lead to an inexact solution (the optimal solution does not satisfy the original equation (3.c)) especially when there is a high reverse current due to the high penetration of RESs.

$$\Re(S_{l,t})^2 + \Im(S_{l,t})^2 \geq f_{l,t} v_{up(l),t} \quad (4)$$

Such an inexact solution could violate the maximum voltage or current limits. In order to ensure the feasibility of the optimization solution, we derive upper bounds on nodal voltage and line current magnitudes and we enforce these uppers bounds to be smaller than the corresponding maximum limits. This idea is inspired by [24], where an exact relaxation method is proposed for the OPF problem, called AROPF. In this paper, we aim to derive a new model for AC power flow that exhibits the same features of AROPF while accounting for the model of the BESS.

$$\underline{S}_{l,t} = s_{l,t} + \underline{p}^b_{l,t} + jq^b_{l,t} + \sum_{m \in \mathcal{L}} G_{l,m} \underline{S}_{m,t} \tag{5.a}$$

$$\bar{v}_{l,t} = \bar{v}_{up(l),t} - 2\Re(z_l^* \underline{S}_{l,t}) \tag{5.b}$$

$$\bar{S}_{l,t} = s_{l,t} + \bar{p}^b_{l,t} + jq^b_{l,t} + \sum_{m \in \mathcal{L}} G_{l,m} \bar{S}_{m,t} + Z_l \bar{f}_{l,t} \tag{5.c}$$

$$\max\{|\Re(\bar{S}_{l,t})|, |\Re(\underline{S}_{l,t})|\}^2 + \max\{|\Im(\bar{S}_{l,t})|, |\Im(\underline{S}_{l,t})|\}^2 \leq \bar{f}_{l,t} v_{up(l),t} \tag{5.d}$$

First, we define auxiliary variable $\bar{v}_{l,t}$ as an upper bound on the variable $v_{l,t}$. For this aim, the square of the current magnitude is assumed to be equal to its lower bound, which is zero. Assuming this, (5.a) and (5.b) express Kirchhoff's current and voltage laws, respectively, where auxiliary variable $\underline{S} = \underline{P} + j\underline{Q}$ represents the corresponding apparent power flow. In this paper, notations without subscript, such as $v$, denote the corresponding matrix, whose $(l,t)^{th}$ element is $v_{l,t}$.

Lemma II: if $(s, p^r, q^b, S, v, f)$ satisfy (3) and $(s, p^r, q^b, \hat{S}, \hat{v})$ satisfy (5.a)-(5.b), then $v \leq \bar{v}$ and $S \geq \underline{S}$.

As the upper bound on the square of the current magnitude, auxiliary variable $\bar{f}$ is introduced in (5.d). Using the Kirchhoff's current law (5.c), we derive also $\bar{S} = \bar{P} + j\bar{Q}$ as an upper bound on the apparent power flow. To define $\bar{f}_{l,t}$, we need to find upper bounds on the absolute value of the active (reactive) power flow. This is equal to $\bar{P}_l(\bar{Q}_l)$ if the active (reactive) power flows in the forward direction (i.e. from $up(l)$ to $l$) and is equal to $-\underline{P}_l(-\underline{Q}_l)$ if the active (reactive) power flows in the reverse direction (see Lemma II). Using these upper bounds on the absolute value of the active and reactive power flows, we derive $\bar{f}$ as in (5.d).

Lemma III: if $(s, p^r, q^b, S, v, f, \underline{S}, \bar{v}, \bar{S}, \bar{f})$ satisfy (3) and (5), then $f \leq \bar{f}$ and $S \leq \bar{S}$.

Lemma II and III are proved in Appendix VII.B. According to these Lemmas, the auxiliary variables $\bar{v}_{l,t}$, $\bar{S}_{l,t}$, and $\bar{f}_{l,t}$ provide upper bounds for the O-OPF variables: $v_{l,t}, S_{l,t}$, and $f_{l,t}$, respectively, and $\underline{S}_{l,t}$ provides a lower bound for $S_{l,t}$.

## IV. Optimization Problem

In this section, the optimization problem is formulated first using the original non-convex models of the battery (1) and grid (3), and then using the developed relaxed convex models. The set of nodes hosting DDGs and BESSs are denoted with $\mathcal{D}$ and $\mathcal{B}$, respectively. The main decision variables are the DDG power setpoints ($s_{i,t} \forall i \in \mathcal{D}$), the reservoir active power setpoint of the BESS ($p^r$), and the reactive power setpoint of the BESS ($q^b$).

### A. Original Optimal Power Flow (O-OPF)

The original formulation of the optimization problem is given in (6), which is referred as Original Optimal Power Flow (O-OPF).

$$F^{obj} = \min_{P_{1,t} \geq 0} \sum_t \mathcal{C}_t \Re(S_{1,t}) \tag{6.a}$$

Subject to:

$$(1.a)-(1.d) \quad : \forall l \in \mathcal{B}, \forall t \tag{6.b}$$

$$(3.a)-(3.c) \quad : \forall l \in \mathcal{L}, \forall t \tag{6.c}$$

$$\begin{cases} v_l^{min} \leq v_{l,t} \\ v_{l,t} \leq v_l^{max} \\ f_{l,t} \leq I_l^{max} \\ \Re(S_{l,t}) \leq P_l^{max} \\ \Im(S_{l,t}) \leq Q_l^{max} \end{cases} \quad : \forall l \in \mathcal{L}, \forall t \tag{6.d}$$

$$\begin{cases} f^b_{l,t} \leq I_l^{b,max} \\ SOC_l^{min} \cdot E_l^{cap} \leq E_{l,t} \\ E_{l,t} \leq SOC_l^{max} \cdot E_l^{cap} \\ E_{l,T}^{DSC} \leq N_l^{DSC} \cdot E_l^{cap} \\ p_{l,dsc}^{min} \leq p^b_{l,t} \leq p_l^{max} \\ \mathcal{R}_i^{dn} \leq p^r_{l,t} - p^r_{l,t-1} \leq \mathcal{R}_i^{up} \end{cases} \quad : \forall l \in \mathcal{B}, \forall t \tag{6.e}$$

$$\begin{cases} p_l^{min} \leq \Re(s_{l,t}) \leq p_l^{max} \\ |s_{l,t}|^2 \leq s_{l,DG}^{max} \end{cases} \quad : \forall l \in \mathcal{D}, \forall t \tag{6.f}$$

where, $\mathcal{C}_t$ is the cost of the energy purchase from the grid during $(t, t+1)$, $v_l^{min}$ and $v_l^{max}$ are the minimum and maximum limits of the squared voltage magnitude, $I_l^{max}$ denotes the maximum limits of the squared current magnitude, $P_l^{max}$ and $Q_l^{max}$ are the maximum active and reactive power limits of line $l$, respectively, $I_l^{b,max}$ is the maximum limit of the squared current magnitude associated with the inverter of the BESS, $E_l^{cap}$ is the energy capacity of the BESS, $SOC_l^{min}$ and $SOC_l^{max}$ represent the minimum and maximum allowable SOC of the BESS, $N_l^{DSC}$ is the maximum allowable number of discharging cycles of the BESS at node $l$ during the planning time horizon, $p_l^{min}$ and $p_l^{max}$ are the minimum and maximum limits of nodal active power consumption, $\mathcal{R}_i^{dn}$ and $\mathcal{R}_i^{up}$ are the maximum allowable rate of power increase and decrease at the reservoir of the BESS, and $s_{l,DG}^{max}$ is the maximum apparent power capacity of the DDG.

As formulated in (6.a), the objective is to minimize the total cost of active power purchase from the grid. Constraints (6.b) and (6.c) refer to the original model of the BESS and grid, respectively. The electrical security constraints including voltage and current limits are accounted for in (6.d). The constraints of the BESS are expressed in (6.e), and finally (6.f) accounts for the modeling of the DDG power capacities. The O-OPF problem is non-convex due to the non-linear constraints (1.b), (1.c), and (3.c). To handle the nonconvexity in the problem, we propose MAROPF approach in the following.

### B. Modified Augmented Relaxed OPF (MAROPF)

In order to convexify O-OPF, the relaxed models of the BESS and grid that were developed, respectively, in section III.B and III.C are applied, resulting in a convex optimization problem, named MAROPF. The main and auxiliary decision variables are $(s, p^r, q^b)$ and $(S, v, f, \bar{S}, \bar{f}, \underline{S}, \bar{v}, \underline{E}, \bar{E}, \bar{E}^c)$, respectively.

$$F^{obj} = \min_{P_{1,t} \geq 0} \sum_t \mathcal{C}_t \Re(S_{1,t}) \tag{7.a}$$

Subject to:

$$(2.\text{a})\text{-}(2.\text{f}) \quad : \forall l \in \mathcal{B}, \forall t \tag{7.b}$$

$$S_{l,t} = s_{l,t} + \bar{p}_{l,t}^b + jq_{l,t}^b + \sum_{m \in \mathcal{L}} G_{l,m} S_{m,t} + Z_l f_{l,t} \tag{7.c}$$

$$(3.\text{b}), (4), (5.\text{a})\text{-}(5.\text{d}) \quad : \forall l \in \mathcal{L}, \forall t \tag{7.d}$$

$$\begin{cases} v_l^{min} \leq v_{l,t} \\ \bar{v}_{l,t} \leq v_l^{max} \\ \bar{f}_{l,t} \leq I_l^{max} \\ \Re(S_{l,t}) \leq \Re(\bar{S}_{l,t}) \leq P_l^{max} \\ \Im(S_{l,t}) \leq \Im(\bar{S}_{l,t}) \leq Q_l^{max} \end{cases} \quad : \forall l \in \mathcal{L}, \forall t \tag{7.e}$$

$$\begin{cases} \bar{f}_{l,t}^b \leq I_l^{b,max} \\ SOC_l^{min} \cdot E_l^{cap} \leq \underline{E}_{l,t} \\ \bar{E}_{l,t} \leq SOC_l^{max} \cdot E_l^{cap} \\ \bar{E}_{l,T}^{DSC} \leq N_l^{DSC} \cdot E_l^{cap} \\ p_l^{min} \leq \hat{p}_{l,t}^b \\ \bar{p}_{l,t}^b \leq p_l^{max} \\ \mathcal{R}_i^{dn} \leq p_{l,t}^r - p_{l,t-1}^r \leq \mathcal{R}_i^{up} \end{cases} \quad : \forall l \in \mathcal{B}, \forall t \tag{7.f}$$

$$\begin{cases} p_{l,DG}^{min} \leq \Re(s_{l,t}) \leq p_{l,DG}^{max} \\ |s_{l,t}|^2 \leq s_{l,DG}^{max} \end{cases} \quad : \forall l \in \mathcal{D}, \forall t \tag{7.g}$$

The relaxed and convex model of the BESS is given in (7.b). The set of constraints (7.c)- (7.d) represents the relaxed model of the grid power flow. It is to be noted that $p_{l,t}^b$ as the non-linear term of (3.a) is replaced with $\bar{p}_{l,t}^b$ in (7.c). The electrical security constraints of the grid are expressed in (7.e), where the maximum voltage and current limits are applied on the upper bounds of the voltage and current magnitudes. Constraints (7.f) represents the set of BESS constraints. The upper and lower bound variables of current, energy, and power are constrained via the corresponding maximum and minimum safety limits.. Finally, the power and energy capacities of the DDGs are formulated in (7.g).

*C. Feasibility Analysis*

In this section, we demonstrate the solution feasibility, i.e., that the exact voltage and current magnitudes respect the grid security constraints (6.d). For the implementation of the proposed MAROPF, the following four conditions should be met, which can be checked a priori. As shown in [24], they are satisfied for most radial distribution networks in practice.

Condition **C1**: $\quad \|\mathbf{E}\| < 1$

Condition **C2**: $\quad \mathbf{DE} \leq \eta_1 \mathbf{D}$

Condition **C3**: $\quad (\mathbf{H}diag(r)\,\mathbf{E}) \circ \mathbf{H} \leq \eta_2\, \mathbf{H}diag(r)$

Condition **C4**: $\quad \mathbf{H}diag(r)\,\mathbf{EE} \leq \eta_3\, \mathbf{H}diag(r)\,\mathbf{E}$

where, $\mathbf{H} = (\mathbf{I} - \mathbf{G})^{-1}$ and $\mathbf{G}$ denotes the adjacency matrix of the directed graph of the distribution network, meaning that $\mathbf{G}_{k,l}$ equals to one if $k = up(l)$, and equals to zero otherwise. The Frobenius norm of a given matrix $\mathbf{A}$ is denoted with $\|\mathbf{A}\|$. $diag(r)$ represents a diagonal matrix, whose $(l,l)$ element is $r_l$. Bold nonitalic letters refer to the matrices that are defined as follows:

$$\mathbf{E} = 2diag(\pi)\mathbf{H}diag(r) + 2diag(\varrho)\mathbf{H}diag(x) + diag(\vartheta)\mathbf{D}$$

$$\mathbf{D} = 2\mathbf{H}diag(r)(\mathbf{H} - \mathbf{I})diag(r) + 2\mathbf{H}diag(x)(\mathbf{H} - \mathbf{I})diag(x) + \mathbf{H}diag(|z|^2)$$

where, $\pi_l = \max\left(\frac{P_l^{max},\ |Hp^{min}|_l}{v_l^{min}}\right)$, $\varrho_l = \max\left(\frac{Q_l^{max},\ |Hq^{min}|_l}{v_l^{min}}\right)$, and $\vartheta_l = \pi_l^2 + \varrho_l^2$.

*Lemma IV*: if $(s, p^r, q^b, S, v, f)$ satisfy O-OPF (6) and $(s, p^r, q^b, S', v', f', \bar{S}', \bar{f}', \underline{S}, \bar{v})$ satisfy MAROPF (7) with $v' \leq v$, then there exists $\bar{f}$ and $\bar{S}$ such that $\bar{f} \leq \bar{f}'$, $\bar{S} \leq \bar{S}'$, and $(s, p^r, q^b, S, v, f, \bar{S}, \bar{f}, \underline{S}, \bar{v})$ satisfies (1)-(5).

The proof of Lemma IV is provided in Appendix VII.C, inspiring from the proof presented in [24]. Lemma IV will be used to prove the solution feasibility, as stated in Theorem I.

*Theorem I*: Under conditions C1-C4, for every feasible solution $(s, p^r, q^b, S, v, f, \bar{S}, \bar{f}, \underline{S}, \bar{v})$ of the MAROPF (7), there exists a feasible solution $(s, p^r, q^b, S^*, v^*, f^*)$ of the O-OPF (6) with the same nodal power injections (s, $p^r, q^b$).

By using Lemma I, Lemma II, and Lemma III, we can prove Theorem I using the same arguments given in [24]. We omit the detailed proof due to space limitations. Theorem I implies that the nodal power injections (s, $p^r, q^b$) that satisfy the MAROPF (7), creates voltage and current magnitudes that are within their feasible regions. Actually, if $(s, p^r, q^b, S, v, f, \bar{S}, \bar{f}, \underline{S}, \bar{v})$ is a feasible solution of (7), then $\bar{p}^b$, $S, f$, and $v$ are not necessarily the exact values of the BESS active power, line power, line current and nodal voltage variables. The exact values can be obtained by feeding the solution of the MAROPF problem $(s, p^r, q^b)$ to the non-linear model of the battery (1) and grid (3).

*D. Optimality Analysis*

*Theorem II*: Under conditions C1-C4, every optimal solution $(s, p^r, q^b, S, v, f, \bar{S}, \bar{f}, \underline{S}, \bar{v})$ of the MAROPF (7) satisfies (2.a) and (4) with strict equality.

According to Theorem II the optimal solution of the MAROPF guarantee the exactness of the relaxed models developed for the battery (2) and grid (4)-(5). It means that if in the formulation of (2.b), $\tilde{p}_{l,t}^r$ is accurately estimated ($\underline{p}_{l,t}^b = p_{l,t}^b$), then the optimal solution of the MAROPF (7) is the global optimal solution of the O-OPF (6). The proof of Theorem II is similar to what is presented in [24], which is omitted here due to space limitations.

In the following, we propose an iterative approach in order to achieve the accurate estimation of $\tilde{p}_{l,t}^r$:

- Step 0: Initialize $\tilde{p}^{r(0)}$ and the iteration index: $k \leftarrow 0$.
- Step 1: Apply the formulation of $\underline{p}^{b(k)}$ as in (2.b).
- Step 2: Solve the MAR-OPF optimization problem (7) and save the solution: $\omega^{(k)} = (s^{(k)}, p^{r(k)}, q^{b(k)})$
- Step3: Update $\tilde{p}^r$ as: $\tilde{p}^{r(k+1)} \leftarrow p^{r(k)}$.
- Step 4: derive the formulation of $p^{b(k)}$ using (1.b). If $\underline{p}^{b(k)} = p^{b(k)}$, stop, and give $\omega^{(k)}$ as the final optimal solution, otherwise increment the iteration index ($k \leftarrow k+1$) and return to Step 1.

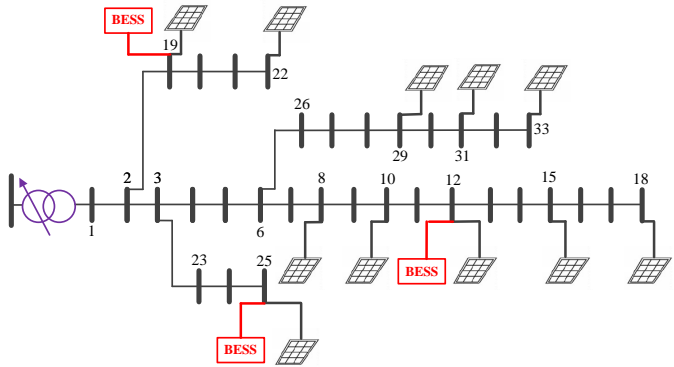

Fig. 4. One line diagram of the modified 34-bus distribution grid [25].

Let assume that the iterative algorithm converges at iteration $k$. The converged solution satisfies $\underline{p}^{b(k)} = p^b$. Moreover, according to Theorem II, at each iteration, the optimal solution satisfies (2.a) and (4) with strict equalities. The only remaining relaxation in the formulation of MAROPF is (2.e), which is related to the maximum limit of the BESS energy. Therefore, if the upper limit of the BESS energy capacity is not binding or if the internal power transmission losses of the BESS is negligible, the obtained global solution upon convergence matches the global optimal solution of O-OPF.

## V. Simulation Results

The developed optimization problem is applied on the IEEE 34-bus standard test distribution network, which is shown in Fig. 4. The nodal and line data are reported in [25]. The test network is modified by placing 11 PV resources at the nodes shown in Fig. 4. Three of these PV resources are equipped with BESSs, whose active and reactive powers are controlled through the control of the interface regulators.

The forecast of the available active power for each PV resource is obtained according to the data related to a sunny day reported in [26]. The load profile at different nodes of the network are according to the industrial, commercial, residential, and rural load patterns reported in [27]. According to ANSI C84.1 standard, the under- and overvoltage limits are set, respectively, to 0.90 and 1.05 p.u. [24].

Fig. 5 shows the nodal voltages during an entire day, when there is no BESS installed in the network. Violation of the voltage limits requires the optimal and secured control of active and reactive powers of each BESS.

All the optimization problems are implemented on a PC with an Intel(R) Xeon(R) CPU and 6 GB RAM; and solved in Matlab/Yalmip environment, using Gurobi solver.

### A. The proposed MAROPF optimization

In this section, the developed MAROPF formulation is applied on the test system. The planning time horizon is assumed from 03:00 A.M until 15:00 P.M with 15 minutes interval. As mentioned in section IV.C, the conditions can be verified a priori to solving the optimization problem. It is illustrated in [29] that these conditions are valid usually for radial distribution networks, including the 34-bus IEEE standard test network.

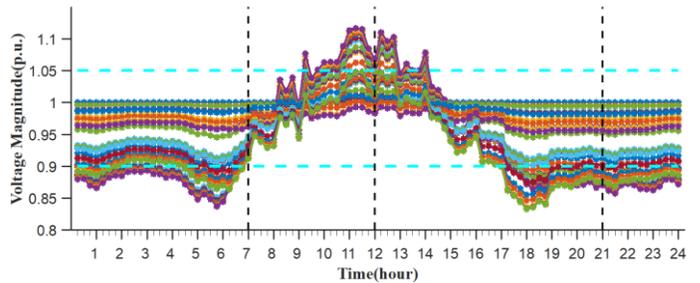

Fig. 5. Voltage profile during the day without any BESS in the 34-bus network. Dashed lines represent voltage limits. The voltages at different buses are shown in different colors.

Table I shows the numerical results obtained from different iterations of the developed MAROPF approach until convergence. As indicated in the table, the proposed algorithm converges with two iterations. As given in the last columns, the proposed algorithm finds the first feasible solution in 2.3 seconds and the best one in 4.7 seconds at the second iteration. This proves that our method can be applied to practical distribution networks with realistic sizes and complexities. The second column of Table I gives the optimal objective values obtained at each iteration. As it can be seen, the quality of the solution is improved through the iterations.

The rest of Table I evaluates the feasibility of the obtained solution at each iteration. For deriving the exact voltage and current magnitudes, the optimal BESS power setpoints are applied on the model of the network in Matlab/MATPOWER toolbox. Then, multi-time power flow simulations are run for all the time steps of the planning horizon. The results of these power flow simulations for voltage and current magnitudes at time 9:30 A.M. (when the voltage reaches its maximum) are illustrated in Fig. 6 and compared with the results of the optimization problem. Although there is a mismatch, the MAROPF guarantees that the exact voltage and current values are within their secured bounds. It should be noted that the objective values reported in the second column of Table I are also corresponding to the exact values that are found with the result of these power flow simulations.

The optimization and power flow simulation results for the state of charge of BESSs during the planning time span are shown in Fig. 7. In addition, the total energy discharge of each BESS is also noted in this figure. For calculating the exact energy level of BESSs, their obtained optimal active and reactive power setpoints $(p^r_{l,t}, q_{l,t})$ are fed into their exact model (1.a)-(1.c). As it can be seen, the lower and upper energy limits and the maximum discharging cycles are all respected. It should be noted that the reported energy values in p.u. are based on 1MWh and not based on each BESS capacity. The mismatch between the exact value of the BESS energy and the value of the variable $\bar{E}_{l,t}$ (labled as BE_opt in Fig. 7) is due to the conservative modelling of BESS ohmic power losses.

Table I. The numerical results of the proposed MAROPF optimization

| Iteration # | $F^{obj}$ (p.u.) | $\min_{i,t}(v_l^{max} - v_{l,t})$ (p.u.) | $\min_{i,t}(v_{l,t} - v_l^{min})$ (p.u.) | $\min_{i,t}(I_l^{max} - f_{l,t})$ (p.u.) | Computation Time (sec) |
|---|---|---|---|---|---|
| 1 | 2.734e-2 | 2.91e-3 | 1.08e-5 | 2.35 | 2.3 |
| 2 | 2.712e-2 | 4.23e-3 | 3.98e-7 | 2.33 | 2.4 |

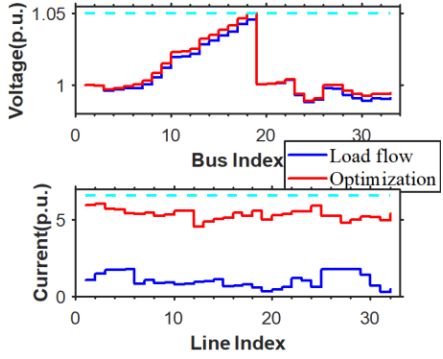

Fig. 6. The results of the MAROPF for the voltage and current profiles at time 9:30 A.M. obtained from the optimization problem and post power flow simulations.

The non-negative margins to all the security limits of the grid and BESSs declare that the solution of the optimization problem at each iteration is feasible. The very small values reported in Table I for some of these margins at the final iteration show that the converged optimal solution is making the most use of the BESS and grid capacity. This could qualitatively suggest that with the iterative MAROPF, the feasible solution space shrinks very little with respect to the O-OPF.

### B. Comparison of alternative models of conversion losses

An claimed in sections IV.CIV.D and IV.D , the proposed approach for modelling conversion losses is robust and accurate, respectively. This has been analytically proved in section VI. In this section, we aim to illustrate this claim by comparison to the alternative approaches in the literature.

#### A.1) Relaxed Linear Model

In this approach, we separately assign to charging and discharging power setpoints of BESS, two nonnegative variables. This approach is applied on the test system with the same conditions as for the proposed MAROPF. The obtained objective value is 2.706e-2, which is lower than the one obtained by MAROP (see Table I). However, ss shown in Fig. 8, the energy of the BESS at node 12 exceeds the upper limit, showing that the obtained solution is actually not feasible. In order to find the reason, the obtained values for the charging and discharging power variables of this BESS are also shown in Fig. 8. As it can be seen, these variables are nonzero at the same time samples. This means that the BESS is simultaneously charging and discharging, which is unrealistic. Sufficient conditions for ensuring non-simultaneous BESS charging and discharging are developed in in [30] and [31] for the case of home energy management and DC-OPF problems, respectively. However, as shown in [30] and [31], these conditions might not be satisfied in some operation conditions.

#### A.2) Binary Model

The approach proposed in [32] addresses the problem of simultaneous non-zero values of charging and discharging power variables by introducing binary variables, indicating if a given BESS is charging or discharging at each time sample. The overall optimization problem is in the form of a Mixed-Integer Second-Order Cone Programming (MISOCP). The integrality constraints in this formulation are handled using the Branch-and-Bound method, with the optimality gap equal to 1e-10 [33].

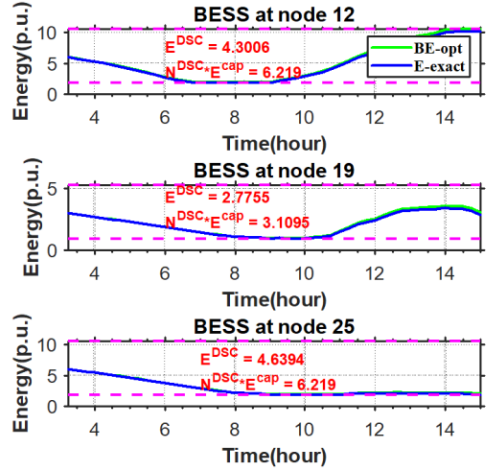

Fig. 7. The MAROPF optimization solution for the state of charge and the discharging energy of BESS.

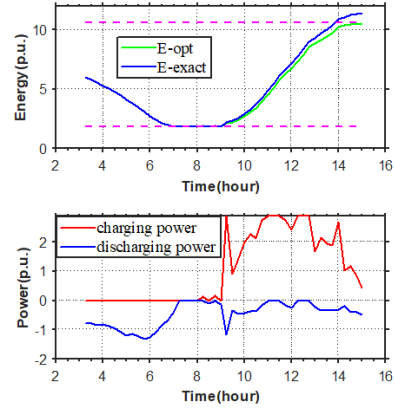

Fig. 8. The power and state of charge of the BESS at bus 12 obtained from the solution of the relaxed linear model for the BESS conversion losses.

The solution of this approach satisfies the BESS energy constraints. However, it is obtained after 59.2 seconds, which is very longer than the proposed MAROP. It is due to the use of binary variables for each BESS and for each time sample. In large scale problems, this computation time could be very longer. The obtained objective value is equal to 2.7245e-2, which is larger than the one obtained from the proposed MAROPF approach. This is due to the relaxation of integrality constraints in the Branch-and-Bound method used for solving the MISOCP problem in the binary approach.

In the literature, Meta-Heuristic methodologies are also applied to model the conversion losses of BESS. Among others are Genetic Algorithm [34], Particle Swarm Optimization [35], and Simulated Annealing [36]. However, there is no guarantee that the obtained solution will be close to the global optimal

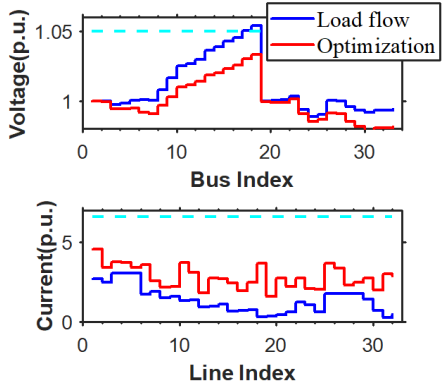

Fig. 9. The results of the ROPF for the voltage and current profiles at time 9:30 A.M. obtained from the optimization problem and post power flow simulations.

solution. It is even possible that an existing feasible solution is not found by these approaches.

*C. Comparison of alternative grid models*

In order to illustrate the performance of the MAROP approach in formulating power flow equations, it is compared with the general relaxation method, named Relaxed Optimal Power Flow (ROPF). This model is constructed by (3.a), (3.b), and (4), where no any auxiliary variables are used for the voltage, power, and currents. Applying this formulation to the test system results in an objective value equal to 1.099e-2, which is less than the one obtained with MAROPF (see Table I). The minimum margin to the upper voltage limit is equal to -3.518, which shows a deviation. Therefore, the obtained solution is not feasible.

For this solution, the exact voltage and current values at time 9:30 A.M. are illustrated in Fig. 9, which are found using post power flow simulations. These profiles are compared with the corresponding values obtained with the ROPF formulation. As shown, there is a mismatch between these profiles, especially at the leaf nodes and lines. This is due to the reverse power flow injected by PVs and BESSs at the leaf nodes that causes the upper voltage limit to bind and, consequently, leading to an inexact solution.

## VI. Conclusion

In this paper, an iterative algorithm is proposed to find the optimal power setpoints of BESSs and DDGs in distribution networks. At each iteration, an SOCP optimization is solved with the objective of minimizing the total cost of power purchase from the grid. For this aim, first, we develop convex relaxation methods for BESSs and grid power flow equations. These methods are then integrated resulting in a convex formulation, named MAROPF. Then, an iterative approach is proposed to reduce the conservativeness of the relaxed BESS model. Analytical proofs are provided to ensure that this iterative algorithm converges with finite number of iterations to a non-conservative solution with regard to the modelling of conversion losses. It is also proved that the solution at each iteration respects all the security constraints of the grid, BESSs and DDGs.

The contributions of this paper were validated on the 32-bus IEEE standard test network. The results of post power flow simulations show that the solution of the proposed algorithm respects all the security constraints of the grid and BESSs. Some of these constraints were met with small margins, which suggests the good quality of the obtained solution. The existing methods for the modelling of BESS conversion losses and for the grid power flow equations are applied on the same test network under the same simulation conditions. The comparison of the obtained results declares the performance of the proposed algorithm in terms of the solution feasibility, quality, and computation time.

## VII. Appendices

In this section, we provide the analytical proof of Lemmas and Theorems presented in sections III and IV. For this aim, first, we rewrite the power flow equations. As shown in [24], the equations (8.a), (8.b), (8.c), and (8.d) are the matrix form of (7.c), (5.a), (3.b) and (5.b), respectively.

$$S = \mathbf{H}(s + \bar{p}^b + jq^b) + \mathbf{H}\text{diag}(Z)f \tag{8.a}$$

$$\underline{S} = \mathbf{H}(s + \underline{p}^b + jq^b) \tag{8.b}$$

$$v = v_0 \mathbf{H}e - 2\mathbf{H}\text{diag}(r)\mathbf{H}(\Re(s) + \bar{p}^b) - 2\mathbf{H}\text{diag}(x)\mathbf{H}(\Im(s) + q^b) - \mathbf{D}f \tag{8.c}$$

$$\bar{v} = v_0 \mathbf{H}e - 2\mathbf{H}\text{diag}(r)\mathbf{H}(\Re(s) + \hat{p}^b) - 2\mathbf{H}\text{diag}(x)\mathbf{H}(\Im(s) + q^b) \tag{8.d}$$

where,

*A. Proof of Lemma I*

If the estimated value of the BESS reservoir power ($\tilde{p}^r$) has the same sign as the exact value ($p^r$), then $\underline{p}^b_{l,t} = p^b_{l,t}$. Otherwise, one of the following two cases may occur.

As the first case, let assume that the BESS is in charging mode. It means that $p^r_{l,t} \geq 0$. Therefore, according to (1.b), $p^b_{l,t} = \frac{p^r_{l,t}}{\eta^c_l}$. Now if $\tilde{p}^r_{l,t} \leq 0$, (1.c) gives $\underline{p}^b_{l,t} = \eta^d_l p^r_{l,t}$. Since $p^r_{l,t} \geq 0$ and $0 \leq \eta^d_l, \eta^c_l \leq 1$, we have $\underline{p}^b \leq p^b$.

In the second case, we consider the discharging state of the BESS, where $p^r_{l,t} \leq 0$. According to (1.b), $p^b_{l,t} = \eta^d_l \cdot p^r_{l,t}$. In this case, if $\tilde{p}^r_{l,t} \geq 0$, we have $\underline{p}^b_{l,t} = \frac{p^r_{l,t}}{\eta^c_l}$. Since $p^r_{l,t} \leq 0$, we conclude that $\underline{p}^b \leq p^b$. Therefore, the result of Lemma I holds in all the possible cases.

*B. Proof of Lemma II and Lemma III*

Lemma II and III are proved in [24] for the case, where the net nodal power vector in the formulation of $(v, f, S)$, $(\hat{S}, \hat{v})$, and $(v, \bar{f}, \bar{S})$ are the same. In our case, the net nodal power in the formulation of $(v, f, S)$ is equal to $s_l + p^b_l + jq^b_l$ (3.a), whereas it is in the formulation of $(\underline{S}, \bar{v})$ is equal to $s_l + \underline{p}^b_l + jq^b_l$ (5.a) and in the formulation of $(v, \bar{f}, \bar{S})$ is equal to $s_l + \bar{p}^b_l + jq^b_l$ (5.c). As proved in [24], with the decrease of the nodal power, $\hat{v}$ is increasing and $\underline{S}$ is decreasing. Since the arguments of Lemma II apply to the case with the same nodal powers, they can also be applied to our case, where $p^b_{l,t}$ in the formulation of $(\underline{S}, \bar{v})$ is replaced with its lower bound $\underline{p}_{l,t}^b$ (Lemma I). In addition, it is proved in [24] that with the increase of the nodal power, $\bar{S}$ is increasing. Since for the case with the same nodal powers, $S \leq \bar{S}$, it can be concluded that this claim applies also to our case. Having this together with the results of Lemma II, we conclude that $\underline{S} \leq S \leq \bar{S}$. Therefore, $|\Re(S)| \leq \max(|\Re(\underline{S})|, |\Re(\bar{S})|)$ and $|\Im(S)| \leq \max(|\Im(\underline{S})|, |\Im(\bar{S})|)$. According to (3.c) and (5.d), we conclude that $f \leq \bar{f}$.

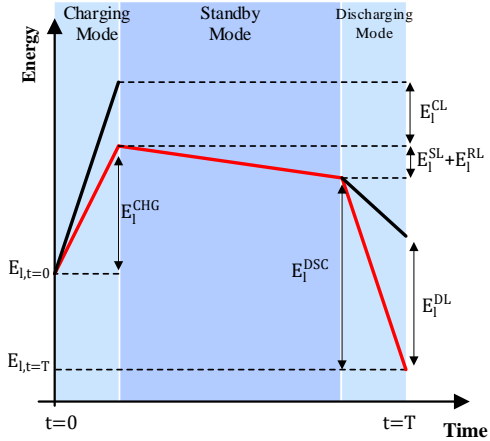

Fig. 10. Energy diagram of BESS throughout the planning time horizon

### C. Proof of Lemma IV

Lemma IV is proved in [24] for the case where in the formulation of $S_l$ and $\bar{S}_l$ the same nodal powers were used ($s_l + p_l^b + jq_l^b$). In that case, it was proved that there exists $\bar{f}$ and $\bar{S}$ such that $\bar{f} \leq \bar{f}'$, $\bar{S} \leq \bar{S}'$. In the developed formulation of MAROPF, the upper-bounds of the nodal powers are sed in (5.c) to derive $\bar{S}_l$. Therefore, according to (5.d), the obtained values of $\bar{f}'$ and $\bar{S}'$ in MAROPF are larger than the ones in the case of [24]. Therefore, the expressions defined in the proof of Lemma IV in [24] for $\bar{f}$ and $\bar{S}$ will be also less than $\bar{f}'$ and $\bar{S}'$ in MAROPF.

### D. Energy Discharge of BESS

In this section, we derive a convex formulation for $E_l^c$ as the total energy discharge of BESS at node $l$ during the whole planning time horizon (from $t = 0$ to $t = T$). The energy diagram of a given BESS through time is depicted in Fig. 10. Without loss of generality, we assume that the BESS is under charging, discharging, and discharging modes only one time through the planning horizon. Moreover, we integrate the ohmic and stand-by losses into the same time period. However, these assumptions do not affect the generality of the following analysis. Equation (9.a) formulates the energy balance, where $E_l^L$ represents the total energy losses of the BESS, which is expressed in (9.b) as the summation of the charging losses ($E_l^{CL}$), discharging losses ($E_l^{DL}$), standby losses ($E_l^{SL}$), and ohmic losses ($E_l^{RL}$). The formulation of $E_l^{CL}$ and $E_l^{DL}$ are given in (9.c) and (9.d), respectively, where $E_l^{CHG}$ refers to the total charging energy of the BESS. Replacing $E_l^{CHG}$ in (9.d) with its equivalent expression in (9.e) and then integrating into (9.a) together with (9.b) and (9.c) leads to (9.f). Moreover, $\Delta E$ can be expressed as in (9.g). Replacing (9.g) in (9.f) leads to (9.h) as the final formulation of $E_l^c$.

$$\Delta E_l = E_{l,t=T} - E_{l,t=0} = \sum_{t=1}^{T} p_{l,t}^b - E_l^L \quad (9.a)$$

$$E_l^L = E_l^{CL} + E_l^{DL} + E_l^{SL} + E_l^{RL} \quad (9.b)$$

$$E_l^{DL} = E_l^{DSC}(1 - \eta^d) \quad (9.c)$$

$$E_l^{CL} = E_l^{CHG}\left(\frac{1}{\eta^c} - 1\right) \quad (9.d)$$

$$E_{l,t=0} + E_l^{CHG} = E_{l,t=T} + E_l^{DSC} + E_l^{SL} + E_l^{RL} \quad (9.e)$$

$$E_l^{DSC} = \left(\frac{\eta^c}{1 - \eta^c\eta^d}\right)\left(\sum_{t=1}^{T} p_{l,t}^b - \frac{E_{SL} + E_l^{Loss} + \Delta E}{\eta^c}\right) \quad (9.f)$$

$$\Delta E_l = \sum_{t=1}^{T} p_{l,t}^r - E_l^{SL} - E_l^{RL} \quad (9.g)$$

$$E_l^{DSC} = \left(\frac{\eta^c}{1 - \eta^c\eta^d}\right)\sum_{t=1}^{T}\left(p_{l,t}^b - \frac{p_{l,t}^r}{\eta^c}\right) \quad (9.h)$$